\documentclass[11pt]{amsart}  
\usepackage{amsmath,amssymb,amsfonts,times,hyperref}

\title{New upper bounds for kissing numbers from semidefinite programming}

\author{Christine Bachoc} 
\address{C. Bachoc, Laboratoire A2X, Universit\'e Bordeaux I, 351,
cours de la Li\-b\'e\-ration, 33405 Talence France}
\email{bachoc@math.u-bordeaux1.fr}

\author{Frank Vallentin} 
\address{F. Vallentin, Centrum voor Wiskunde en Informatica (CWI),
Kruislaan 413, 1098 SJ Amsterdam, The Netherlands}
\email{f.vallentin@cwi.nl}

\thanks{The second author was supported by the Netherlands
Organization for Scientific Research under grant NWO 639.032.203 and
by the Deutsche Forschungsgemeinschaft (DFG) under grant SCHU
1503/4-1.}


\date{October 3, 2007}
\newtheorem{defi}{Definition}[section]

\newtheorem{proposition}[defi]{Proposition}
\newtheorem{theorem}[defi]{Theorem}
\newtheorem{remark}[defi]{Remark}
\newtheorem{corollary}[defi]{Corollary}

\newcommand{\R}{{\mathbb{R}}} 

\newcommand{\F}{{\mathbb{F}}} 
\newcommand{\N}{{\mathbb{N}}}

\newcommand{\I}{{\mathcal{I}}} 
\newcommand{\mymid}{:} 

\newcommand{\Sn}{S^{n-1}} 
\newcommand{\Snm}{S^{n-2}} 
\newcommand{\Pd}{\operatorname{Pol}_{\leq d}}
\newcommand{\Pk}{\operatorname{Pol}_{\leq k}}
\newcommand{\Pki}{\operatorname{Pol}_{\leq k+i}}
\newcommand{\On}{{\operatorname{O}(\R^n)}}
\newcommand{\Onn}{{\operatorname{O}(\R^{n-1})}}
\newcommand{\Omp}{{\operatorname{O}(\R^{m+1})}}

\newcommand{\Trace}{\operatorname{Trace}}
\newcommand{\Harm}{\operatorname{Harm}}
\newcommand{\1}{\operatorname{\bf 1}}

\newcommand{\card}{\operatorname{card}}

\newcommand{\prodeucl}[2]{#1 \cdot #2}
\newcommand{\prodhaar}[2]{(#1,#2)}

\begin{document}

\begin{abstract} 
Recently A. Schrijver derived new upper bounds for binary codes using
semidefinite programming. In this paper we adapt this approach to
codes on the unit sphere and we compute new upper bounds for the
kissing number in several dimensions. In particular our computations
give the (known) values for the cases $n = 3, 4, 8, 24$.
\end{abstract}

\maketitle

\section{Introduction}

In geometry, the kissing number problem asks for the maximum number
$\tau_n$ of unit spheres that can simultaneously touch the unit sphere
in $n$-dimensional Euclidean space without pairwise overlapping.  The
value of $\tau_n$ is only known for $n=1,2,3,4,8,24$. While its
determination for $n=1,2$ is trivial, it is not the case for other
values of $n$. 

The case $n=3$ was the object of a famous discussion between Isaac
Newton and David Gregory in 1694. For a historical perspective of this
discussion we refer to \cite{C}. The first valid proof of the fact
``$\tau_3=12$'', like in the icosahedron configuration, was only given
in 1953 by K.~Sch\"utte and B.L.~van der Waerden in \cite{SW}.

In the seventies, P.~Delsarte developed a method, initially aimed at
bounding codes on finite fields (see \cite{D}), that yields an upper
bound for $\tau_n$ as a solution of a linear program and more
generally yields an upper bound for the size of spherical codes of
given minimal distance. We shall refer to this method as the LP
method. With this method, A.M.~Odlyzko and N.J.A.~Sloane (\cite{OS}),
and independently V.I.~Levenshtein (\cite{Le}), proved $\tau_8=240$
and $\tau_{24}=196560$ which are respectively the number of shortest
vectors in the root lattice $E_8$ and in the Leech lattice. For other
values of $n$, the LP method gives in many cases the best known upper
bounds.  However, for $n=3$ and $n=4$ it only gives the upper bounds
$\tau_3\leq 13$ and $\tau_4\leq 25$.

In 2003, O.R.~Musin succeeded to prove the conjectured value
$\tau_4=24$, which is the number of shortest vectors in the root
lattice $D_4$, with a variation of the LP method (see \cite{M} and the
survey \cite{PZ} of F.~Pfender and G.M.~Ziegler).

To complete the picture, let us discuss uniqueness of the optimal
point configurations. For dimensions $8$ and $24$, uniqueness was
proved by E.~Bannai and N.J.A.~Sloane (\cite{BS}). Their proof
exploits the fact that the LP method obtains exactly the aimed value.
For dimension $3$, there are infinitely many possible configurations.
In the regular icosahedron configuration, the angular distances
between the contact points are strictly greater than the required
$\pi/3$, hence these points can be moved around obtaining infinitely
many new suitable configuration. This partially explains why the
determination of $\tau_3$ is difficult. On the contrary, uniqueness of
the optimal configuration of points in dimension $4$ is widely
believed, but remains unproven.

The LP method, which was established by P.~Delsarte, J.M.~Goethals and
J.J.~Seidel in \cite{DGS}, handles the more general problem of the
determination of a bound for the maximal number
\[
A(n,\theta) = \max\{\card(C) \mymid \mbox{$C \subset \Sn$ with
$\prodeucl{c}{c'} \leq \cos \theta$ for $c, c' \in C$, $c \neq c'$}\}
\]
of points on the unit sphere with minimal angular distance $\theta$.
Such configurations of points, also called {\em spherical codes with
minimal angular distance $\theta$}, are of special interest in
information theory. The kissing number problem is equivalent to the
problem of finding $A(n,\pi/3)$.

In this paper, we define a semidefinite program (SDP for short) whose
optimal solution gives an upper bound for $A(n,\theta)$ and
strengthens the LP method. Computational results show that for several
values of $n$ this SDP method gives better upper bounds for $\tau_n$
than the LP method.

To be more precise, let us recall that the LP method relies on the
existence of polynomials $P_k^n(t)$, satisfying the so-called
positivity property:
\begin{equation}\label{pos 1}
\text{for all finite }C\subset \Sn,
\sum_{(c,c')\in C^2} P_k^n(\prodeucl{c}{c'})\geq 0.
\end{equation}
These polynomials arise as zonal spherical polynomials on the sphere,
i.e.\ the zonal polynomials associated to the decomposition of the
space of polynomial functions under the action of the orthogonal group
$\On$.

The consideration of the action restricted to a subgroup $H$ of
$\On$, chosen to be the stabilizer group of a fixed point $e\in
\Sn$, leads us to some symmetric matrices $S_k^n$ whose coefficients
are symmetric polynomials in three variables such that
\begin{equation}\label{pos 2}
\text{ for all finite }C\subset \Sn, \sum_{(c,c',c'')\in C^3}
S_k^n(\prodeucl{c}{c'}, \prodeucl{c}{c''}, \prodeucl{c'}{c''})\succeq
0
\end{equation}
where the sign ``$\succeq 0$'' stands for: ``is positive
semidefinite''. The reason why we obtain matrices instead of functions
comes from the fact that, in the decomposition of the space of
polynomial functions on the sphere under the action of $H$,
multiplicities greater than $1$ appear. In fact these multiplicities
are exactly the sizes of the corresponding matrices.  From \eqref{pos
1} and \eqref{pos 2} we derive an SDP whose solution gives an upper
bound for $A(n,\theta)$.

Our approach adapts the method proposed by A.~Schrijver in \cite{S} to
the unit sphere whereas he obtains new upper bounds for binary codes
from an SDP. His work can also be interpreted in group theoretic
terms, involving the isometry group of the Hamming space $\F_2^n$ and
the subgroup stabilizing $(0,\dots,0)$ which is the group of
permutations of the $n$ positions. It is very likely that many other
spaces of interest in coding theory can be treated likewise. The case
of non-binary codes was considered by D.C.~Gijswijt, A.~Schrijver and
H.~Tanaka in \cite{DST}.

The paper is organized as follows: Section \ref{LP} reviews on the LP
method.  Section \ref{SD} introduces and calculates the semidefinite
zonal matrices associated to the action of $H$ and leading to the
matrices $S_k^n$.  Section \ref{SDP} defines the semidefinite program
and its dual that establishes the desired bound.  Section
\ref{results} discusses computational results.

\section{Review of the LP Method on the Unit Sphere}
\label{LP}

We introduce the following notations. The standard inner product of
the Euclidean space $\R^n$ is denoted by $\prodeucl{x}{y}$. The unit
sphere
\[
\Sn := \{x\in \R^n \mymid \prodeucl{x}{x}=1\}
\]
is homogeneous under the action of the orthogonal group $\On = \{O \in
\R^{n \times n} : O^t O = I_n\}$, where $I_n$ denotes the identity
matrix. It is moreover two-point homogeneous, meaning that the orbits
of $\On$ on pairs of points are characterized by the value of their
inner product. The space of real polynomial functions of degree at
most $d$ on $\Sn$ is denoted by $\Pd(\Sn)$. It is endowed with the
induced action of $\On$, and equipped with the standard
$\On$-invariant inner product
\[
\prodhaar{f}{g}=\frac{1}{\omega_n}\int_{\Sn} f(x)g(x)d\omega_n(x),
\] 
where $\omega_n = \frac{2\,\pi^{n/2}}{\Gamma(n/2)}$ is the surface area of $\Sn$
for the standard measure $d\omega_n$.  It is a well-known fact (see
e.g. \cite[Ch.~9.2]{VK}) that under the action of $\On$
\begin{equation}\label{dec 1}
\Pd(\Sn)=H_0^n\perp H_1^n\perp\ldots\perp H_d^n,
\end{equation}
where $H_k^n$ is isomorphic to the $\On$-irreducible space
\[
\Harm_k^n = \Big\{f \in \R[x_1, \ldots, x_n] \mymid \text{$f$
homogeneous}, \deg f = k, \sum_{i=1}^n \frac{\partial^2}{\partial
x_i^2} f = 0\Big\}
\]
of harmonic polynomials in $n$ variables which are homogeneous and
have degree $k$.  We set $h_k^n:=\dim(\Harm_k^n) = \binom{n+k-1}{n-1}
- \binom{n+k-3}{n-1}$.

A certain family of orthogonal polynomials is associated to the unit
sphere. They will be denoted by $P_k^n$, with the convention that
$P_k^n$ has degree $k$ and is normalized by $P_k^n(1)=1$. For $n \geq
3$ these polynomials are up to multiplicative constants Gegenbauer
polynomials $C^{\lambda}_k$ with parameter $\lambda = n/2-1$.  So they
are given by $P_k^n(t) = C^{n/2-1}_k(t)/C^{n/2-1}_k(1)$, and the
Gegenbauer polynomials $C^{\lambda}_k$ can be inductively defined by
$C_0^\lambda(t) = 1$, $C_1^\lambda(t) = 2\lambda t$, and
\[
k C_k^\lambda(t) = 2(k + \lambda - 1) t C^{\lambda}_{k-1}(t) -
(k+2\lambda-2) C^{\lambda}_{k-2}(t), \quad \text{for $k \geq 2$}.
\]
They are orthogonal with respect to the weight function
$(1-t^2)^{\lambda-1/2}$ on the interval $[-1,1]$.  For $n = 2$ the
polynomials $P_k^n$ coincide with the Chebyshev polynomials of the
first kind $T_k$ which can be inductively defined by $T_0(t) = 1$,
$T_1(t) = t$, and
\[
T_k(t) = -2t T_{k-1}(t) + T_{k-2}(t), \quad \text{for $k \geq 2$},
\]
and they are orthogonal with respect to the weight function
$(1-t^2)^{-1/2}$ on the interval $[-1,1]$.

The polynomials $P_k^n(t)$ are related to the decomposition (\ref{dec
1}) by the so-called {\em addition formula} (see
e.g. \cite[Ch.~9.6]{AAR}): for any orthonormal basis
$(e_1,\ldots,e_{h_k^n})$ of $H_k^n$ and for any pair of points $x,y
\in \Sn$ we have
\begin{equation}\label{add form 1}
P_k^n(\prodeucl{x}{y})=\frac{1}{h_k^n}\sum_{i=1}^{h_k^n}e_i(x)e_i(y).
\end{equation}
From the addition formula (\ref{add form 1}), the positivity property
(\ref{pos 1}) becomes obvious:
\begin{align*} 
\sum_{(c,c')\in C^2}P_k^n(\prodeucl{c}{c'})&=\sum_{(c,c')\in C^2}\frac{1}{h_k^n}\sum_{i=1}^{h_k^n}e_i(c)e_i(c')\\
&=\frac{1}{h_k^n}\sum_{i=1}^{h_k^n}\sum_{(c,c')\in
  C^2}e_i(c)e_i(c')=\frac{1}{h_k^n}\sum_{i=1}^{h_k^n}\left(\sum_{c\in C}e_i(c)\right)^2\geq 0.
\end{align*}

Now we introduce the unknowns of the LP to be considered. For a spherical code $C$ 
we define the two point distance distribution
\begin{equation*}\label{xu}
x(u):=\frac{1}{\card(C)}\card\{(c,c')\in C^2\mymid \prodeucl{c}{c'}=u\},
\end{equation*}
where $u \in [-1,1]$. Clearly, only a finite number of $x(u)$'s are not
equal to zero, and the positivity property can be rewritten as a
linear inequality in the $x(u)$'s:
\begin{equation}\label{pos P}
\sum_{u \in [-1,1]} x(u) P_k^n(u)\geq 0.
\end{equation}
Moreover, the number of elements of $C$ is given by
$\card(C)=\sum_{u \in [-1,1]} x(u)$. Noticing the obvious conditions
$x(1)=1$, $x(u)\geq 0$, and $x(u)=0$ for $\cos\theta<u<1$ if
the minimal angular distance of $C$ is $\theta$, we are led to
consider the following linear program: For any $d \geq 1$, the optimal
solution of the linear program
\begin{equation}\label{LP primal}
\begin{array}{ll}
\max\big\{ & \displaystyle 1 + \sum_{u\in [-1,\cos\theta]}x(u) \quad \mymid \\[2ex]
& \quad\qquad \mbox{$x(u)=0$ for all but finitely many $u\in [-1,\cos\theta]$,}\\[1ex]
& \quad\qquad \mbox{$x(u) \geq 0$ for all $u\in [-1,\cos\theta]$,}\\[1ex]
& \quad\qquad \mbox{$1+\sum_{u\in [-1,\cos\theta]} x(u) P_k^n(u) \geq
0$ for all $k=1,\ldots, d$}
\big\},
\end{array}
\end{equation}
gives an upper bound for $A(n,\theta)$.  The dual linear problem is
\begin{equation}\label{LP dual}
\begin{array}{ll}
\min\big\{ & 1 + \displaystyle \sum_{k=1}^d f_k \quad \mymid \\
& \quad\qquad \mbox{$f_k \geq 0$ for all $k = 1, \ldots, d$},\\[1ex]
& \quad\qquad \mbox{$\sum_{k=1}^d f_k P_k^n(u) \leq -1$ for all $u\in [-1,\cos\theta]$}\big\}.
\end{array} 
\end{equation}
By the duality theorem (cf.~\cite{Du}) any feasible
solution of (\ref{LP dual}) gives an upper bound for the optimal
solution of (\ref{LP primal}).  The dual linear program can be
restated in the following way involving polynomials:

\begin{theorem}\label{Th LP} (see e.g.\ \cite[Th.~4.3]{DGS}, \cite{KL}, \cite{OS}, \cite[Ch.~9]{CS})

\noindent Let $F(t)=\sum_{k=0}^d f_k P_k^n(t)$ be a polynomial of
degree at most $d$ in $\R[t]$. If 

(a) $f_k\geq 0$\  for all $k\geq 1$ and $f_0>0$ and

(b) $F(u)\leq 0$\  for all $u\in [-1,\cos\theta]$,

\noindent then
\[\displaystyle A(n,\theta)\leq \frac{F(1)}{f_0}.\]
\end{theorem}

\section{Semidefinite Zonal Matrices}
\label{SD}

Now we fix a point $e\in \Sn$, and let $H:=\mathrm{Stab}(\On,e)$
be the stabilizer of $e$ in $\On$. Obviously, $H\simeq \Onn$ since $\Onn$
can be identified with the orthogonal group of the orthogonal
complement of $\R e$. 

It is a classical result (see e.g. \cite[Ch.~9.2]{VK}) that for the
restricted action to $H$ the decomposition of $\Harm_k^n$ into
$H$-irreducible subspaces is given by:
\[
\Harm_k^n\simeq \bigoplus_{i=0}^k \Harm_i^{n-1}.
\]
Hence, each of the $H_k^n$ in \eqref{dec 1} decomposes likewise:
\begin{equation}
\label{decomp hkn}
H_k^n=H_{0,k}^{n-1}\perp H_{1,k}^{n-1}\perp\ldots \perp H_{k,k}^{n-1},
\end{equation}
where $H_{i,k}^{n-1}\simeq \Harm_i^{n-1}$. We give an explicit
description of this decomposition in the proof of Theorem~\ref{Th Y}.

We summarize the situation in the following picture.

\[
\begin{array}{ccccccccc}
\Pd(\Sn) & = & H_0^n         & \perp & H_1^n         & \perp & \ldots & \perp & H_d^n\\
         & = & H_{0,0}^{n-1} & \perp & H_{0,1}^{n-1} & \perp & \ldots & \perp & H_{0,d}^{n-1}\\
         &   &               & \perp & H_{1,1}^{n-1} & \perp & \ldots & \perp & H_{1,d}^{n-1}\\
&   &               &       &               & \multicolumn{4}{r}{\cdots\cdots\cdots\cdots\cdots\cdots}\\
         &   &               &       &               &        &  & \perp & H_{d,d}^{n-1}\\
\end{array}
\]

\noindent The isotypic components of the $H$-decomposition of $\Pd(\Sn)$ are 
\begin{equation}\label{dec Ik}
\I_k:=H_{k,k}^{n-1} \perp \ldots \perp H_{k,d}^{n-1}\simeq
(d-k+1)\Harm_k^{n-1}, \quad \mbox{for $k = 0, \ldots, d$.}
\end{equation}
Now we show how to associate to each $\I_k$ a ``zonal matrix'' in view
of an analogue of the addition formula (\ref{add form 1}).

\begin{theorem}\label{Th Z}
Let $\I=R_0\perp R_1\perp\ldots\perp R_m\simeq (m+1)R$ be an isotypic
component of $\Pd(\Sn)$ under the action of $H$, with $R$ an
$H$-irreducible space of dimension $h$. Let $(e_{0,1},\ldots,e_{0,h})$
be an orthonormal basis of $R_0$ and let $\phi_i:R_0\to R_i$ be
$H$-isomorphisms preserving the inner product on $\Pd(\Sn)$. Let
$e_{i,j}=\phi_i(e_{0,j})$, so that $(e_{i,1},\ldots,e_{i,h})$ is an
orthonormal basis of $R_i$. Define
\[E(x):=\Big(\frac{1}{\sqrt{h}}e_{i,j}(x)\Big)_{\substack{0\leq i\leq m\\1\leq j\leq h}}=
\frac{1}{\sqrt{h}}\begin{pmatrix}
e_{0,1}(x) & \ldots & e_{0,h}(x)\\
\vdots&&\vdots\\
e_{m,1}(x) & \ldots & e_{m,h}(x)
\end{pmatrix},
\]
and
\[Z(x,y):=E(x)E(y)^t \in \R^{(m+1) \times (m+1)}.\]
Then the following properties hold for the matrix $Z$:

(a) $Z(x,y)$ does not depend on the choice of the orthonormal basis
  of $R_0$.

(b) The change of $\phi_i$ to $-\phi_i$ for some $i$ or the choice
  of another decomposition of $\I$ as a sum of $m+1$ orthogonal
  $H$-submodules changes
  $Z(x,y)$ to some $OZ(x,y)O^t$ with $O\in \Omp$.

(c) For all $g\in H$, $Z(g(x),g(y))=Z(x,y)$.

(d) (Matrix-type positivity property) 
\begin{equation}\label{Z pos}
\text{For all finite }C\subset \Sn, \sum_{(c,c')\in C^2}Z(c,c')\succeq
0.
\end{equation}
\end{theorem}

\begin{proof}
(a) If $(\epsilon_{0,1}, \ldots,\epsilon_{0,h})$ is another
orthonormal basis of $R_0$, then there is an orthogonal $h \times h$
matrix $O$ with $(\epsilon_{0,1}, \ldots,\epsilon_{0,h})=(e_{0,1},
\ldots,e_{0,h})O$.  In this case the matrix $E(x)$ is changed to
$E(x)O$ and, since $OO^t=I_h$, the matrix $Z(x,y)$ stays unchanged.

(b) By Schur's Lemma and by the irreducibility of $R$, there are only
two possible choices for $\phi_i$, namely $\phi_i$ and $-\phi_i$, once
the subspaces $R_i$ are fixed.

Let $\I=S_0\perp\ldots\perp S_m$ be another decomposition of $\I$,
together with $H$-isomorphisms $\psi_i:R_i\to S_i$ preserving the
inner product on $\Pd(\Sn)$. Then $\psi=(\psi_0,\ldots,\psi_m)$
defines an $H$-endomorphism of $\I$. Again by Schur's Lemma, for a
suitable choice of basis in $R_i$ and by permuting rows and columns,
the matrix of $\psi$ is block diagonal with $h$ blocks of size $(m+1)
\times (m+1)$ and with the same $(m+1)\times (m+1)$ matrix $O\in
\Omp$ as blocks. This means that $E(x)$ changes to $OE(x)$ and
so $Z(x,y)$ becomes $OZ(x,y)O^t$.

(c) Since $e_{i,j}(g^{-1}(x))=(g e_{i,j})(x)$, the computation of
$Z(g^{-1}(x),g^{-1}(y))$ amounts to replace in the definition of
$Z(x,y)$ the $e_{i,j}$ by $g e_{i,j}$.  Since $R_i$ is $H$-stable,
$\epsilon_{i,j}:=g e_{i,j}$, with $j = 1, \ldots, h$, is another
orthonormal basis of $R_i$, and
\[
\phi_i(\epsilon_{1,j})=\phi_i(g
   e_{1,j})=g\phi_i(e_{1,j})=g e_{i,j}=\epsilon_{i,j}.
\]
Hence from (a) we conclude $Z(g^{-1}(x),g^{-1}(y))=Z(x,y)$.

(d) We have $\displaystyle \sum_{(c,c')\in C^2}Z(c,c')=\Big(\sum_{c\in
C}E(c)\Big)\Big(\sum_{c\in C}E(c)\Big)^t\succeq 0 $.
\end{proof}

The orbits of $H$ on pairs of points on the unit sphere $x,y \in \Sn$
are characterized by the values of the three inner products
$\prodeucl{e}{x}$, $\prodeucl{e}{y}$, and $\prodeucl{x}{y}$. By
definition the coefficients $Z_{i,j}(x,y)$ of $Z(x,y)$ are polynomials
in the variables $x_1, \ldots, x_n, y_1, \ldots, y_n$. Then, property
(c) of Theorem~\ref{Th Z} implies that they can be expressed as
polynomials in the three variables $u = \prodeucl{e}{x}$, $v =
\prodeucl{e}{y}$, and $t = \prodeucl{x}{y}$.

By $Z_k^n$, for $0\leq k\leq d$, let us denote the matrix associated
to $\I_k$ as defined above, and more precisely to the decomposition
\eqref{dec Ik} of $\I_k$. Now we shall calculate the matrix
$Y_k^n(u,v,t)$ with
\begin{equation}\label{add form 2}
Z_k^n(x,y)=Y_k^n(\prodeucl{e}{x},\prodeucl{e}{y}, \prodeucl{x}{y}).
\end{equation}

\begin{theorem}\label{Th Y}
With the above notations, we have, for all $0\leq i,j\leq d-k$,
\begin{equation}\label{Y}
\big(Y_k^n\big)_{i,j}(u,v,t)= \lambda_{i,j} P_{i}^{n+2k}(u)P_{j}^{n+2k}(v)Q_k^{n-1}(u,v,t),
\end{equation}
where 
\[
\displaystyle Q_k^{n-1}(u,v,t):=\big((1-u^2)(1-v^2)\big)^{k/2}P_k^{n-1}\Big(\frac{t-uv}{\sqrt{(1-u^2)(1-v^2)}}\Big),
\]
and
\[
\lambda_{i,j}=\frac{\omega_n}{\omega_{n-1}}\frac{\omega_{n+2k-1}}{\omega_{n+2k}}(h_i^{n+2k}h_j^{n+2k})^{1/2}.
\]
\end{theorem}

\begin{proof}
We explicitly use an orthonormal basis of $H_{k,k+i}^{n-1}$ to calculate
$Y_k^n(u,v,t)$. Such a basis is constructed in
\cite[Ch.~9.8]{AAR}. Let us recall the construction. For $x\in \Sn$,
let
\[
x=ue+\sqrt{1-u^2}\zeta,
\]
where $u=\prodeucl{x}{e}$ and $\zeta$ belongs to the unit sphere
$\Snm$ of $(\R e)^{\perp}$. With $f\in H_k^{n-1} \subset \Pk(\Snm)$ we
associate $\varphi(f)\in \Pk(\Sn)$ defined by:
\[
\varphi(f)(x)=(1-u^2)^{k/2}f(\zeta).
\]
Note that the multiplication by $(1-u^2)^{k/2}$ forces $\varphi(f)$ to
be a polynomial function in the coordinates of $x$. Clearly $\varphi$
commutes with the action of $H$. Hence $\varphi(H_k^{n-1})$ is a
subspace of $\Pk(\Sn)$ which is isomorphic to $\Harm_k^{n-1}$. More
generally, the set $\{\varphi(f)P(u) : f \in \Harm_{k}^{n-1}, \deg P
\leq i\}$ is a subspace of $\Pki(\Sn)$ which is isomorphic to $i+1$
copies of $\Harm_{k}^{n-1}$.  By induction on $k$ and $i$ there exist
polynomials $P_i(u)$ of degree $i$ such that
$\varphi(H_k^{n-1})P_i(u)=H_{k,k+i}^{n-1}$. Note that this
construction could be used to derive decomposition \eqref{decomp hkn}
explicitly.

We can exploit the fact that the subspaces $H_{k,l}^{n-1}$ are
pairwise orthogonal to prove an orthogonality relation between the
polynomials $P_i$. Then this orthogonality relation will enable us to
identify the polynomials $P_i$ as multiples of Gegenbauer polynomials.
Let us recall that the measures on $\Sn$ and on $\Snm$ are related by:
\[
d\omega_n(x)=(1-u^2)^{(n-3)/2}du d\omega_{n-1}(\zeta).
\]
Whenever $i \neq j$ we have for all $f\in H_k^{n-1}$
\begin{align*}
0&=\frac{1}{\omega_n}\int_{\Sn} \varphi(f)P_i(u)\varphi(f)P_j(u) d\omega_n(x)\\
&=\frac{1}{\omega_n}\int_{\Sn} f(\zeta)^2(1-u^2)^kP_i(u)P_j(u) d\omega_n(x)\\
&=\frac{1}{\omega_n}\int_{\Snm} f(\zeta)^2d\omega_{n-1}(\zeta)\int_{-1}^1 (1-u^2)^{k+(n-3)/2}P_i(u)P_j(u) du,
\end{align*}
from which we derive that 
\[
\int_{-1}^1 (1-u^2)^{k+(n-3)/2}P_i(u)P_j(u) du=0;
\]
hence the polynomials $P_i(u)$ are proportional to $P_i^{n+2k}(u)$.
We obtain an orthonormal basis of $H_{k,k+i}^{n-1}$ from an
orthonormal basis $(f_1,\ldots,f_h)$ of $H_k^{n-1}$ by taking
$e_{i,j}=\lambda_i\varphi(f_j)P_i^{n+2k}(u)$ for a suitable normalizing
factor $\lambda_i$. We compute $\lambda_i$ in a similar way as above:
\begin{align*}
1&=\frac{1}{\omega_n}\int_{\Sn} \big(\lambda_i\varphi(f_j)P_i^{n+2k}(u)\big)^2 d\omega_n(x)\\
&=\frac{1}{\omega_n}\int_{\Snm} \big(f_j(\zeta)\big)^2d\omega_{n-1}(\zeta)\int_{-1}^1
\lambda_i^2(1-u^2)^{k+(n-3)/2}\big(P_i^{n+2k}(u)\big)^2 du\\
&=\frac{\omega_{n-1}}{\omega_n}\int_{-1}^1 \lambda_i^2(1-u^2)^{k+(n-3)/2}\big(P_i^{n+2k}(u)\big)^2 du.
\end{align*}
From the addition formula \eqref{add form 1} applied to
$(P_i^{n+2k}(u))^2$ one easily shows that
\[
\int_{-1}^1
(1-u^2)^{k+(n-3)/2}\big(P_i^{n+2k}(u)\big)^2du=
\frac{\omega_{n+2k}}{\omega_{n+2k-1} h_{i}^{n+2k}},
\]
so we obtain
\[
\lambda_i^2=\frac{\omega_n}{\omega_{n-1}}\frac{\omega_{n+2k-1}}{\omega_{n+2k}}h_{i}^{n+2k}.
\]
Now we are in the situation of Theorem \ref{Th Z} with 
\[
R_0=H_{k,k}^{n-1}, R_1 = H_{k,k+1}^{n-1}, \ldots, R_{d-k}=H_{k,d}^{n-1}
\]
and their orthonormal basis $(e_{0,1},\ldots,
e_{0,h}),\ldots,(e_{d-k,1},\ldots, e_{d-k,h})$. The isomorphisms
$\phi_i$ are the multiplications by
$(\lambda_i/\lambda_1)P_i^{n+2k}(u)$.

Then, the coefficient $(i,j)$, with $0\leq i,j\leq d-k$, of $Z_k^n$ is
given by: 

\begin{eqnarray*}
& & \big(Z_k^n\big)_{i,j}(x,y)\\
&=&\frac{1}{h}\sum_{s=1}^h e_{i,s}(x)e_{j,s}(y)\\
&=&\frac{1}{h}\sum_{s=1}^h  \lambda_i(1-u^2)^{k/2}f_s(\zeta)P_i^{n+2k}(u)\lambda_j(1-v^2)^{k/2}f_s(\xi)P_j^{n+2k}(v)\\
&=&\lambda_i\lambda_j P_i^{n+2k}(u)P_j^{n+2k}(v)\big((1-u^2)(1-v^2)\big)^{k/2}\frac{1}{h}\sum_{s=1}^h f_s(\zeta)f_s(\xi)\\
&=& \lambda_i\lambda_j P_i^{n+2k}(u)P_j^{n+2k}(v)\big((1-u^2)(1-v^2)\big)^{k/2}P_k^{n-1}(\prodeucl{\zeta}{\xi}),
\end{eqnarray*}
where we have written $y=ve+\sqrt{1-v^2}\xi$ and where we applied the
addition formula \eqref{add form 1} to get the last equality. Now we
define $\lambda_{i,j} = \lambda_i \lambda_j$ and since
\[
\prodeucl{\zeta}{\xi}=(t-uv
)/\sqrt{(1-u^2)(1-v^2)},
\]
we have completed the proof.
\end{proof}

\begin{remark} We would like to point out that the role of the number
  $d$ is only to cut $Y_k^n$ to a matrix of finite size. Indeed, $d$
  does not enter in the expression of $\big(Y_k^n\big)_{i,j}(u,v,t)$.
  It is better to view the matrices $Y_k^n$ as matrices of infinite
  size with all finite principal minors having the matrix-type
  positivity property.
\end{remark}

\begin{remark}
  For the semidefinite programming bounds in Section~\ref{SDP} we only
  use the matrix-type positivity property of the matrices $Y^n_k$.
  This property is preserved if one replaces $Y^n_k$ by $A Y^n_k A^t$
  with an invertible matrix $A$. So, e.g.\ one could replace the
  expression of $\big(Y_k^n\big)_{i,j}(u,v,t)$ in \eqref{Y} by the
  simpler $u^i v^j Q_k^{n-1}(u,v,t)$.
\end{remark}

Due to the specific choice of the unit vector $e$ defining the
subgroup $H$, the coefficients of $Y_k^n$ are not symmetric
polynomials. We introduce the symmetrization $S_k^n$ of $Y_k^n$ and
state the announced property \eqref{pos 2}.

\begin{corollary} For all $d\geq 0$, for all $k\geq 0$, let $Y_k^n$ be
  the matrix in Theorem \ref{Th Y} and let $S_k^n$ be defined by
\begin{equation}
S_k^n=\frac{1}{6} \sum_{\sigma} \sigma Y_k^n,
\end{equation}
where $\sigma$ runs through the group of all permutations of the
variables $u,v,t$ which acts on matrix coefficients in the obvious
way.  Then the matrices $S_k^n$ are symmetric and have symmetric
polynomials as coefficients.  We have:
\begin{equation}\label{pos Y}
\text{For all finite $C\subset \Sn$,} \sum_{(c,c')\in
  C^2}Y_k^n(\prodeucl{e}{c},\prodeucl{e}{c'}, \prodeucl{c}{c'})\succeq 0,
\end{equation}
and
\begin{equation}\label{pos S}
\text{for all finite $C\subset \Sn$,} \sum_{(c,c',c'')\in
  C^3}S_k^n(\prodeucl{c}{c'}, \prodeucl{c}{c''},\prodeucl{c'}{c''})\succeq 0.
\end{equation}
\end{corollary}

\begin{proof}
Note that $\big(Y_k^n\big)_{j,i}(u,v,t)=\big(Y_k^n\big)_{i,j}(v,u,t)$
which gives the desired properties of $S_k^n$. Property \eqref{Z pos}
rephrases to \eqref{pos Y} and property \eqref{pos S} is obtained from
\eqref{pos Y} by taking $e=c''\in C$ and summing over all $c''\in C$.
\end{proof}

To end this section we show that the positivity property \eqref{pos 1}
is actually a consequence of the matrix-type positivity property
\eqref{pos 2}. As shown in the following proposition one can express
the polynomials $P_k^n$ as a linear combination of diagonal elements
of the matrices $Y_k^n$ with non-negative coefficients.

\begin{proposition}
\label{PY}
We have the following expression for the polynomials $P_k^n(t)$ in
terms of matrix coefficients of $Y_k^n(u,v,t)$:
\begin{equation}
P_k^n(t)=\sum_{s=0}^k \frac{h_s^{n-1}}{h_k^n}\big(Y_s^n\big)_{k-s,k-s}(u,v,t).
\end{equation}
Consequently, property \eqref{pos Y} or property \eqref{pos 2} implies \eqref{pos 1}.
\end{proposition}

\begin{proof}
The addition formula \eqref{add form 1} holds for any orthonormal
basis of $H_k^n$.  We take an orthonormal basis of $H_k^n$ obtained by
concatenation of orthormal basis of the spaces $H_{0,k}^{n-1}$,
$H_{1,k}^{n-1}, \dots, H_{k,k}^{n-1}$. If $(e_{s,1}^k,
e_{s,2}^k,\dots, e_{s,h_{s}^{n-1}}^k)$ denotes an orthonormal basis of
$H_{s,k}^{n-1}$, we have from Theorem \ref{Th Z}
\begin{equation*}
\big(Y_s^n\big)_{k-s,k-s}(e\cdot x, e\cdot y,x\cdot y)=\frac{1}{h_s^{n-1}}\sum_{i=1}^{h_s^{n-1}}e_{s,i}^k(x)e_{s,i}^k(y).
\end{equation*}
By the addition formula \eqref{add form 1}
\begin{align*}
P_k^n(x\cdot y) &= \frac{1}{h_k^n}
\sum_{s=0}^{k}\sum_{i=1}^{h_s^{n-1}}e_{s,i}^k(x)e_{s,i}^k(y)\\
&=\sum_{s=0}^{k} \frac{h_s^{n-1}}{h_k^n}\big(Y_s^n\big)_{k-s,k-s}(e\cdot x, e\cdot y,x\cdot y),
\end{align*}
and hence
\begin{equation*}
P_k^n(t)=\sum_{s=0}^k \frac{h_s^{n-1}}{h_k^n}\big(Y_s^n\big)_{k-s,k-s}(u,v,t).
\end{equation*}

Since the coefficients ${h_s^{n-1}}/{h_k^n}$ are non-negative, and
since the diagonal elements of a semidefinite matrix are non-negative,
\eqref{pos 1} is a consequence of \eqref{pos Y}.

With the action of the permutation group of the variables $u,v,t$
\begin{equation*}
\frac{1}{3}\big(P_k^n(u)+P_k^n(v)+P_k^n(t)\big)=\sum_{s=0}^k \frac{h_s^{n-1}}{h_k^n}\big(S_s^n\big)_{k-s,k-s}(u,v,t).
\end{equation*}
Replacing $u=c\cdot c'$, $v=c\cdot c''$, $t=c'\cdot c''$ and summing
over $(c,c',c'')\in C^3$ for a code $C$, we obtain \eqref{pos 1} from
\eqref{pos 2}.

\end{proof}

\section{The Semidefinite Programming Bound}\label{SDP}

In this section we set up an SDP whose optimum gives an upper bound
for $A(n,\theta)$ which is at least as good as the LP method.

For a spherical code $C$ we define the three points distance
distribution
\[
x(u,v,t):=\frac{1}{\card(C)}\card\{(c,c',c'')\in C^3\mymid
\prodeucl{c}{c'}=u, \prodeucl{c}{c''}=v, \prodeucl{c'}{c''}=t\},
\]
where $u,v,t\in [-1,1]$ and the matrix 
\[
\begin{pmatrix} 1&u&v\\u&1&t\\v&t&1\end{pmatrix},
\]
being the Gram matrix of three vectors on a unit sphere, is positive
semidefinite.

The last condition together with the first is equivalent to the fact
that the determinant of the Gram matrix is non-negative, hence
\begin{equation}\label{det}
1+2uvt-u^2-v^2-t^2 \geq 0.
\end{equation}
The two point distance distribution $x(u)$ as defined in Section
\ref{LP} and the three point distance distribution $x(u,v,t)$ are
related by $x(u,u,1) = x(u)$. The three point distance distribution
satisfies the following obvious properties:
\[
\begin{array}{ll}
& \displaystyle x(u,v,t)\geq 0,\\[1ex]
& \displaystyle x(1,1,1)=1,\\[1ex]
& \displaystyle
x(\sigma(u),\sigma(v),\sigma(t))=x(u,v,t)\;\;\mbox{for
  all permutations $\sigma$,}\\[1ex]
& \displaystyle\sum_{u,v,t} x(u,v,t)=\card(C)^2,\\
& \displaystyle\sum_{u} x(u,u,1)=\card(C).\\
\end{array}
\]
Furthermore, from the positivity properties \eqref{pos P} and
\eqref{pos S}, we have for any $d \geq 0$:
\[
\begin{array}{ll}
& \displaystyle \sum_{u} x(u,u,1)P_k^n(u)\geq 0\;\;\mbox{for
  $k = 1, \ldots, d$,}\\
& \displaystyle \sum_{u,v,t} x(u,v,t)S_k^n(u,v,t)\succeq
0\;\;\mbox{for $k = 0, \ldots, d$,}
\end{array}
\]
where the matrix $S_k^n$ has size $(d-k+1) \times (d-k+1)$.
If the minimal angular distance of $C$ is $\theta$, we have moreover
\[
x(u,v,t)=0\;\;\mbox{whenever $u,v,t\notin [-1,\cos \theta] \cup \{1\}$.}
\]
To factor out the action of the permutations of the variables $u,v,t$
we introduce the domains
\[
D=\{(u,v,t)\mymid \mbox{$-1\leq u\leq v\leq t\leq \cos\theta$ and $1+2uvt-u^2-v^2-t^2\geq 0$\},}\]
\[D_0=\{(u,u,1)\mymid -1\leq u\leq \cos\theta \},\qquad 
I=[-1,\cos\theta],
\]
and $m(u,v,t)$ with
\[
m(u,v,t)=\begin{cases}
\;\mbox{$6\;$ if $u\neq v\neq t$,}\\
\;\mbox{$3\;$ if $ u=v\neq t$ or $u\neq v=t$ or $u=t\neq v$,}\\
\;\mbox{$1\;$ if $u=v=t$.}
\end{cases}
\]
From the discussion above, a solution to the following semidefinite program in the variables
$x'(u,v,t)=m(u,v,t)x(u,v,t)$
is an upper bound for $A(n,\theta)$:
\[
\begin{array}{ll}
1 + & \max\big\{ \displaystyle \frac{1}{3}\sum_{u\in I} x'(u,u,1)
\mymid\\
& \mbox{$\displaystyle x'(u,v,t)=0\;$  for all but finitely many $(u,v,t)\in D\cup D_0$,}\\[1ex]
& \mbox{$\displaystyle x'(u,v,t)\geq 0\;$ for all $(u,v,t)\in D\cup D_0$,}\\[1ex]
& \displaystyle \left(\begin{smallmatrix} 1&0\\0&0 \end{smallmatrix}\right)+ 
\frac{1}{3}\sum_{u\in I} x'(u,u,1)\left(\begin{smallmatrix} 0&1\\1&1
\end{smallmatrix}\right)+
\sum_{(u,v,t)\in D} x'(u,v,t)\left(\begin{smallmatrix} 0&0\\0&1
\end{smallmatrix}\right)\succeq 0, \\
& \mbox{$\displaystyle 3+\sum_{u\in I} x'(u,u,1)P_k^n(u)\geq 0\;$ for $k=1,\ldots,d$,}\\
& \mbox{$\displaystyle S_k^n(1,1,1)+\sum_{(u,v,t)\in D\cup D_0}
  x'(u,v,t)S_k^n(u,v,t)\succeq 0\;$ for $k=0,\ldots,d$} \big\}.
\end{array}
\]
The third constraint deserves some further explanation. We have
already noticed that
\[
\card(C)^2=1+\sum_{(u,v,t)\in D\cup D_0} x'(u,v,t)=\Big(1+\sum_{u\in
    I}x(u,u,1)\Big)^2,
\]
which implies
\[
\sum_{(u,v,t)\in D} x'(u,v,t)+\frac{1}{3}\sum_{u\in
    I}x'(u,u,1)-\Big(\frac{1}{3}\sum_{u\in I}x'(u,u,1)\Big)^2\geq 0,
\]
and this is equivalent to the semidefinite condition:
\[
\begin{pmatrix}
\displaystyle 1 & \displaystyle \frac{1}{3}\sum_{u\in I}x'(u,u,1)\\
\displaystyle \frac{1}{3}\sum_{u\in I}x'(u,u,1) & \displaystyle \sum_{(u,v,t)\in D}
x'(u,v,t)+\frac{1}{3}\sum_{u\in I}x'(u,u,1)
\end{pmatrix} \succeq 0.
\]

\begin{remark} We want to point out that, despite of the fact that
  \eqref{pos 2} implies \eqref{pos 1}, as is proved in Proposition
  \ref{PY}, the inequalities $3+\sum_{u\in I} x'(u,u,1)P_k^n(u)\geq
  0\;$ should not be removed from our SDP.  Indeed, the last
  inequalities do not imply them for an arbitrary set of numbers
  $x'(u,v,t)$, unless these numbers satisfy the additional equalities:
\begin{equation*}
\sum_{u,v} x(u,v,t) =\big(\sum_{u} x(u,u,1) \big) x(t,t,1) \quad\text{ for
  all }t.
\end{equation*}
These equalities do hold for codes, but they are not semidefinite
conditions. It can be noticed that the third constraint in the
maximization problem above is a weaker consequence of them.
\end{remark}

Just like in the LP method, the main problem with the above SDP, is
that the unknowns $x(u,v,t)$ are indexed by a continuous domain of
$\R^3$. We cannot exploit the information that only a finite number of
them are not equal to zero, because we don't know to which values of
$(u,v,t)$ they correspond. We solve this problem by applying duality
theory. 

Before we derive the SDP dual to the above one we recall the principle
of weak duality. We use the standard notation for the inner product of
symmetric matrices: $\langle A, B\rangle=\Trace(AB)$. Let $J$ be a
(possible infinite) set of indices, let $S_j \in \R^{m \times m}$ be
symmetric matrices with $j \in J$, let $C \in \R^{m \times m}$ be a
symmetric matrix, and let $c_j \in \R$ be real numbers. Suppose that
the real numbers $x_j \in \R$ are a feasible solution of the primal
problem, i.e.\ $x_j = 0$ for all
but finitely many $j \in J$, and $C - \sum_{j \in J} x_j S_j \succeq
0$. Furthermore, suppose that the symmetric matrix $F \in \R^{m \times
  m}$ is a feasible solution of the dual problem, i.e.\ $\langle F,
S_j\rangle = c_j$ for all $j \in J$, and $F \succeq 0$. Then, we have
$\sum_{j \in J} c_j x_j = \langle \sum_{j \in J} x_j S_j, F \rangle
\leq \langle C, F \rangle.  $

In our case this specializes as follows: The set of indices is $J =
D_0 \cup D$. The matrices $S_{(u,v,t)}$ are block matrices with four
blocks of different type. We get one block for each positivity
constraint in the above SDP. So $F$ is also a block matrix with four
blocks of different type. In this case it can be simplified to three
blocks. The first block of $F$ consists of the matrix
$\left(\begin{smallmatrix} b_{11} & b_{12}\\ b_{12} & b_{22}
  \end{smallmatrix}\right)$. The second block of $F$ is the diagonal
matrix with coefficients $a_1, \ldots, a_d$ The third block of $F$ is
again a block matrix with blocks $F_0, \ldots, F_d$ which have the
same size as the matrices $S^n_k$. The matrix $C$ is a block matrix as
well. The first block of $C$ contains the matrix
$\left(\begin{smallmatrix} 1 & 0\\0 & 0\end{smallmatrix}\right)$. The
first entry of the second block is $3$, the other entries in this
block are zero. The third block of $C$ consists of the matrices
$S^n_k(1,1,1)$.  The real numbers $c_{(u,v,t)}$ are zero if $(u,v,t)
\in D$ and equal to $1/3$ if $(u,v,t) \in D_0$. In the following
theorem we give the SDP dual to the above one.  Furthermore we apply
the simplification $S_k^n(1,1,1)=0$ for $k\geq 1$.

\begin{theorem}\label{dual SDP}
Any feasible solution of the following semidefinite problem gives an
upper bound on $A(n,\theta)$:
\[
\begin{array}{ll}
1 + \min\big\{ & \displaystyle \sum_{k=1}^d a_k  + b_{11} + \langle
F_0, S_0^n(1,1,1) \rangle \mymid \\
& \left(\begin{smallmatrix} b_{11} & b_{12}\\ b_{12} & b_{22} \end{smallmatrix}\right) \succeq 0,\\
& \mbox{$a_k \geq 0\;$ for $k = 1,\ldots,d$,}\\
& \mbox{$F_k \succeq 0\;$ for $k = 0,\ldots,d$,}\\
& \displaystyle \sum_{k=1}^d a_k P_k^n(u) + 2b_{12} +
b_{22} +3 \sum_{k=0}^d \langle F_k, S_k^n(u,u,1)\rangle\leq -1,\\
& \displaystyle b_{22} + \sum_{k=0}^d \langle F_k, S_k^n(u, v,t)\rangle \leq 0\big\},
\end{array} 
\]
where the last inequality holds for all $(u,v,t)\in D$ and the second
to last inequality holds for all $u \in I$.
\end{theorem}

Note that if the last inequality holds for all $(u,v,t) \in D$, then
it also holds for the larger domain
\[
D' := \{(u,v,t) \mymid \mbox{$-1\leq u, v, t\leq \cos\theta$ and
$1+2uvt-u^2-v^2-t^2\geq 0$}\},
\]
because the coefficients in $S_k^n$ are symmetric polynomials.

\section{Computational Results}
\label{results}

In this section, we describe one possible strategy to derive explicit
upper bounds for $\tau_n$ from Theorem \ref{dual SDP}. Thereby we make
use of techniques form polynomial optimization introduced
e.g.\ in \cite{L} and \cite{Pa} which we shall briefly recall here.

We consider the polynomials
\[
p(u) = -(u+1/4)^2 + 9/16,
\]
\[
p_1(u,v,t) = p(u),\quad p_2(u,v,t) = p(v),\quad p_3(u,v,t) = p(t),
\]
\[
p_4(u,v,t) = 1 + 2uvt - u^2 - v^2 - t^2,
\]
and we obviously have
\begin{eqnarray*}
I & = & \{u \in \R : p(u) \geq 0\},\\
D' & = & \{(u,v,t) \in \R^3 : p_i(u,v,t) \geq 0,\; i = 1, \ldots, 4\}.
\end{eqnarray*}

We say that a polynomial $f \in \R[x_1, \ldots, x_n]$ is a {\em sum of
squares} if it can be written as $f = \sum_{i=1}^k g_i^2$, for $k \in
\N$ and $g_i \in \R[x_1, \ldots, x_n]$. A polynomial $p(x_1, \ldots,
x_n)$ of degree $2m$ is a sum of squares if and only if there is a
positive semidefinite matrix $Q$ so that $p(x_1, \ldots, x_n) = z^t Q
z$ where $z$ is the vector of monomials $z = (1, x_1, \ldots, x_n, x_1
x_2, \ldots, x_{n-1} x_n, \ldots, x_n^m)$. So assuring that a
polynomial is a sum of squares is a semidefinite condition.

It is easy to see that the last two conditions of the semidefinite
program in Theorem~\ref{dual SDP} are satisfied if the following two
equalities hold:
\[
\begin{array}{l}
- 1 - \displaystyle \sum_{k=1}^d a_k P_k^n(u) - 2b_{12} -
b_{22} -3 \sum_{k=0}^d \langle F_k, S_k^n(u,u,1)\rangle = q(u) + p(u) q_1(u),\\
-\displaystyle b_{22} - \sum_{k=0}^d \langle F_k, S_k^n(u, v,t)\rangle = r(u, v, t) + \sum_{i = 1}^4 p_i(u,v,t) r_i(u,v,t)
\end{array}
\]
where $q, q_1$ and $r, r_1, \ldots, r_4$ are sums of squares. 

It is not apriori clear that the relaxation of using this specific
sum of squares representation is strong enough. The following theorem
of M.~Putinar justifies our approach.

\begin{theorem} (\cite{Pu})
\label{putinar}
Let $K = \{x \in \R^n : p_1(x) \geq 0, \ldots, p_s(x) \geq 0\}$ be a
compact semialgebraic set. Suppose that there is a polynomial $P$ of
the form $P = q + p_1 q_1 + \cdots + p_s q_s$, where $q$ and all
$q_i$'s are sums of squares, so that the set $\{x \in \R^n: P(x) \geq
0\}$ is compact. Then, every polynomial $p$ which is positive on $K$
can be written as $p = r + p_1 r_1 + \cdots + p_s r_s$, where $r$ and
all $r_i$'s are sums of squares.
\end{theorem}

Now we use these considerations to formulate a finite-dimensional
semidefinite program which gives an upper bound on the kissing number
$\tau_n$: We fix $d$ and restrict the polynomials $q, q_1, r, r_1,
\ldots, r_4$ to polynomials having degree at most $N$, with $N \geq
d$. Then we can use the computer to find a feasible solution of this
finite-dimensional semidefinite program.  A feasible solution of it is
at the time a feasible solution of the SDP in Theorem~\ref{dual SDP}.
So it gives an upper bound on the kissing number $\tau_n$.

We implemented this approach and give our results in Table~5.1.

\begin{table}[htb]
\label{Table kissing}
\begin{tabular}{c|c|l|c|c}
    & best lower  & best upper bound                       & LP     & SDP\\
$n$ & bound known & previously known & method & method \\
\hline
3  & 12  & \ 12\hspace{1.7ex} (Sch\"utte, v.d.\ Waerden \cite{SW})  & 13  & 12 \\
4  & 24  & \ 24\hspace{1.7ex} (Musin \cite{M}) & 25  & 24 \\
5  & 40  & \ 46\hspace{1.7ex} (Odlyzko, Sloane \cite{OS}) & 46  & 45 \\
6  & 72  & \ 82\hspace{1.7ex} (O., S. \cite{OS}) & 82  & 78 \\
7  & 126 & 140\hspace{1ex}  (O., S. \cite{OS}) & 140 & 135 \\
8  & 240 & 240\hspace{1ex}  (O., S. \cite{OS}, Levenshtein \cite{L}) & 240 & 240 \\
9  & 306 & 379\hspace{1ex}  (Rzhevskii, Vsemirnov \cite{RV}) & 380 & 366 \\ 
10 & 500 & 594\hspace{1ex} (Pfender \cite{Pf}) & 595 & 567 \\
\end{tabular}
\\[0.3cm]
Table 5.1. Bounds on $\tau_n$.
\end{table}

The values of the last column were found by solving the above
semidefinite program for the values $d = 10$ and $N = 10$. The values
of the third column were obtained by Odlyzko and Sloane by
Theorem~\ref{Th LP} using the value $d = 30$. They pointed out that
even $d = 11$ would suffice for $n \leq 10$. Our calculations were
performed by the program {\tt csdp} (Version 5.0) of B.~Borchers
(\cite{B}) which is available on the Internet
(\href{http://infohost.nmt.edu/~borchers/csdp.html}{http://infohost.nmt.edu/\~{}borchers/csdp.html}).
After solving the SDP with {\tt csdp} we checked independently whether
the solution satisfies the desired constraints. This can be done using
rational arithmetic only. So our computations give rigorous proofs of
the stated upper bounds. Due to numerical instabilities we were not
able to perform this calculation for larger $n$ and/or larger $d$,
$N$. The smallest values of $d$ and $N$ which solve the kissing number
problem in dimension~$3$ is $d = N = 5$. Then, we obtain by the SDP
method $\tau_3 \leq 12.8721$. For the kissing number problem in
dimension~$4$ it is $d = N = 7$, and the SDP method gives $\tau_4 \leq
24.5797$.

For the lower bounds in the first column we refer to the Catalogue of
Lattices of G.\ Nebe and N.J.A.\ Sloane
(\href{http://www.research.att.com/~njas/lattices/kiss.html}{http://www.research.att.com/\~{}njas/lattices/kiss.html}).

Using the polynomial $p(u) = -(u+1/3)^2 + 4/9$ we computed upper
bounds for $A(n, \cos^{-1} 1/3)$. Hereby we improved several entries
of the Table 9.2 of \cite{CS} where all best upper bounds previously
known were obtained by the LP method. We give our results in
Table~5.2. Again we used the values $d = 10$ and $N = 10$ to obtain
the last column.

\begin{table}[htb]
\begin{tabular}{c|c|c|c}
    & best lower  & best upper bound & SDP \\
$n$ & bound known & previously known & method \\
\hline
3 & 9 & 9 & 9 \\
4 & 14 & 15 & 15 \\
5 & 20 & 24 & 23 \\
6 & 32 & 37 & 35 \\
7 & 56 & 56 & 56 \\
8 & 64 & 78 & 74 \\
9 & 96 & 107 & 99 \\ 
10 &  & 146 & 135 \\
\end{tabular}
\\[0.3cm]
Table 5.2. Bounds on $A(n,\cos^{-1} 1/3)$.
\end{table}

We were also able to improve the best known upper bounds for the
so-called Tammes problem with $N$ spheres: What is the largest minimal
angle $\theta(N)$ that can be obtained by a spherical code of $S^2$
with cardinality $N$. Let us recall that the answer is only known for
$N\leq 12$ and for $N=24$ (see \cite[Ch.~1]{CS}). For $N=13$, the
best known lower bound is $0.997223593\approx 57.1367031^{\circ}$
whereas the best-known upper bound is $1.02746114\lessapprox
58.8691870^{\circ}$ due to K.~B\"or\"oczky and L.~Szabo \cite{BSz}. We
obtained $A(3, \cos^{-1}(0.5225)) \leq 12.99$ using $d = N = 10$,
giving the new upper bound of $1.02101593\lessapprox
58.4999037^{\circ}$.  Other values are collected in Table~5.3; the
lower bounds are taken from the homepage of N.J.A.~Sloane
(\href{http://www.research.att.com/~njas/packings}{http://www.research.att.com/\~{}njas/packings/}). The
upper bounds for $N\geq 14$ where established in \cite{BSz2}.

\begin{table}[htb]
\begin{tabular}{c|c|c|c}
    & best lower  & best upper bound & SDP \\
$N$ & bound known & previously known & method \\
\hline
13 & 57.13 & 58.87 & 58.50\\
14 & 55.67 & 58.00 & 56.58\\
15 & 53.65 & 55.84 & 55.03\\
16 & 52.24 & 53.92 & 53.27\\
17 & 51.09 & 52.11 & 51.69\\
\end{tabular}
\\[0.3cm]
Table 5.3. Bounds on $\theta(N)$ (given in degrees).
\end{table}

\section*{Acknowledgements}

We thank Eiichi Bannai, Tatsuro Ito, Monique Laurent, Oleg Musin,
Florian Pfender, Lex Schrijver, and Achill Sch\"urmann for valuable
discussions on this topic.  We thank the anonymous referee for helpful
comments and suggestions.

\end{document}